\documentclass[english,a4paper,reqno]{amsart}

\usepackage{a4}
\usepackage{babel}
\usepackage{fontenc}
\usepackage{amsfonts}
\usepackage[latin1]{inputenc}
\usepackage{amssymb}
\usepackage{amsmath}
\usepackage[bookmarks=true]{hyperref}
\newcommand{\R}{\mathbb{R}}
\newcommand{\N}{\mathbb{N}}

\newcommand{\C}{\mathbb{C}}

\newcommand{\Z}{\mathbb{Z}}

\newcommand{\Li}{\mbox{Li}}
\newtheorem{definition}{Definition}[section]
\newenvironment{Proof}{\noindent {\bf{Proof.}}}{$\Box$\\}
\newtheorem{proposition}{Proposition}[section]

\numberwithin{equation}{section}

\begin{document}

\title[Recurrence relations for the Lerch $\Phi$ function and applications]{Recurrence relations for the Lerch $\Phi$ function and
applications}

\author{Marco Dalai}

%\maketitle
\address{Department of Electronics for Automation, University of Brescia, Via Branze 38, 25123 Brescia (Italy)}
\email{marco.dalai@ing.unibs.it}



\begin{abstract}
In this paper we present a simple method for deriving
recurrence relations and we apply it to obtain two equations
involving the Lerch Phi function and sums of Bernoulli and Euler
polynomials. Connections between these results and those obtained
by H.M. Srivastava, M.L. Glasser and V. Adamchik (\cite{SriGlaAdam2000})
are pointed out, emphasizing the usefulness of this approach with
some meaningful examples.
\end{abstract}
\maketitle

\section{Introduction}

In recent years the properties of the Riemann Zeta function for
positive integer values of its argument have received a lot of
attention. Rapidly convergent series representation such as
Euler's well-known one

\begin{equation}\label{EulSer}
  \zeta(3)=-\frac{4 \pi^2}{7}\sum_{k=0}^\infty
  \frac{\zeta(2k)}{(2k+1)(2k+2)2^{2k}}
\end{equation}
\\
have been recently discovered by several authors in many different
ways (see \cite{Sri_Taiwan2000} for an exhaustive overview).

In this contest, H.M Srivastava, M.L. Glasser e V. S. Adamchick
(\cite{SriGlaAdam2000}) found some interesting series representations of
the values $\zeta(2n+1)$, $n\in\N$, studying
different possible evaluations of the definite integral

\begin{eqnarray}
  I_s(\omega) & = & \label{def:Iom1} \int_0^{\pi/\omega} t^s\csc^2t\;dt\qquad\\
              & = & \label{def:Iom2} -\left(\frac{\pi}{\omega}\right)^s
                    \cot\left(\frac{\pi}{\omega}\right)+ \int_0^{\pi/\omega} t^{s-1}\cot
                    t\;dt\\
              &   & \nonumber\\
              &   & (\Re(s)>1;\quad\omega>1)\nonumber
.
\end{eqnarray}
At the end of their article, the authors demonstrate that every
representation they found can be derived from the following
unification  formula
\begin{align}\label{EquSriv}
  \sum_{k=0}^\infty \frac{\zeta(2k)}{(2k+p)\omega^{2k}}& =\frac{\pi
                  i}{2(p+1)\omega}-\frac{1}{2}\log\left(1-e^{2\pi i
                  /\omega}\right)\\\nonumber
                & -\frac{p\,!}{2}\left(\frac{i\omega}{2 \pi}\right)^p
                  \zeta(p+1)+\frac{1}{2}\sum_{k=1}^p \binom{p}{k}k!\left(\frac{i  \omega}{2
                  \pi}\right)^k\mbox{Li}_k\left( e^{2 \pi
                  i/\omega}\right)\\
                &  \qquad(p\in \N: \quad \omega\in\R,\,
                |\omega|>1)\nonumber.
\end{align}

This paper presents a very simple method for deriving some
recurrence relations and shows how it can be used as a \textquotedblleft short cut{\textquotedblright}
 to
obtain two formulas that generalize equation (\ref{EquSriv}).
The usefulness of the obtained results is then shown by deducting in a simple and unified way many known (or equivalent to known) evaluations and series representations for the Riemann Zeta function and the Dirichlet Beta function.

\section{Notations}
This section presents the notation used in the next sections.
\\The Euler Gamma function $\Gamma(s)$ is defined as the analytic
continuation of the integral
\begin{equation}\label{def:gamma}
   \Gamma(s)=\int_0^\infty  t^{s-1}e^{-t}dt\qquad (s>1)
\end{equation}
and, for $n\in\N$, satisfies the property $\Gamma(n+1)=n!$
\\The Hurwitz Zeta function $\zeta(s,b)$ is defined as the analytic
continuation of the following series (defined for $\Re (s) >1$ and
$\Re(b)>0$)
\begin{equation}\label{def:zeta}
  \zeta(s,b) = \sum_{k=0}^{\infty}\frac{1}{(k+b)^s}
\end{equation}
and it reduces to the Riemann Zeta function $\zeta(s)$ in the case
$b=1$. Similarly, the Dirichlet Beta function $\beta(s)$ is the
analytic continuation of the series
\begin{equation}\label{def:beta}
  \beta(s) = \sum_{k=0}^{\infty}\frac{(-1)^k}{(2k+1)^s}\qquad (\Re(s)>0).
\end{equation}
For $a\in\C$ with $|a|\le1$, taken $s\in\C$ satisfying $\Re(s)>1$
if $a=1$ and $ \Re(s)>0$ if $|a|=1 \wedge a\neq1$, the
Polylogarithm function is defined by
\begin{gather}\label{def:Li}
\Li_s(a) =  \sum_{k=1}^{\infty}\frac{a^k}{k^s}
\end{gather}
and the Lerch $\Phi$ function is defined, for $\Re(b)>0$, by
\begin{gather}\label{def:LerchPhi}
  \Phi(a,s,b)  =  \sum_{k=0}^{\infty}\frac{a^k}{(k+b)^s}\,.
\end{gather}
Furthermore we use Bernoulli polynomials $B_k(x)$ defined by
\begin{equation}\label{def:PolBern}
  \frac{te^{xt}}{e^t-1}=\sum_{k=0}^{\infty}B_k(x)\frac{t^k}{k!}\qquad
  (|t|<2\pi)
\end{equation}
and Euler polynomials $E_k(x)$ defined by
\begin{equation}\label{def:PolEul}
  \frac{2e^{xt}}{e^t+1}=\sum_{k=0}^{\infty}E_k(x)\frac{t^k}{k!}\qquad
  (|t|<\pi).
\end{equation}
The $B_k(0)$ values, called Bernoulli numbers, are represented as
$B_k$, while Euler numbers are the $2^k E_k(1/2)$ values,
represented as $E_k$.
\\Two important (and well known) results that relate Bernoulli and Euler numbers
with the values $\zeta(2k)$ and $\beta(2k+1)$, where $k\in\N$, are
the following relations (\cite{handbook})
\begin{equation}\label{ident:zetanumbern}
 \zeta(2k)=\frac{(-1)^{k-1}(2\pi)^{2k}B_{2k}}{2(2k)!}\,,
\end{equation}
\bigskip
\begin{equation}\label{ident:betanumeul}
 \beta(2k+1)=\frac{(-1)^{k}(\pi/2)^{2k+1}E_{2k}}{2(2k)!}\,.
\end{equation}
Finally, the basic notions of complex variable analysis will be
used.

\section{Holomorphic functions in a strip}

\begin{proposition}\label{teo:funol}Let $\alpha$ and $\beta$ be positive
real numbers and let $f(z)$ be an holomorphic function in the
strip $S= \left\{z\in\C\;:\;\Re(z)\geq
0,\Im(z)\in(-\alpha,\beta)\right\}$ such that $f(z)=O(z^\nu)$,
where $\nu<0$, if $\Re(z)\to\infty$. Then the following equality
holds
\begin{equation}\label{ident:funol}
  \int_0^\infty f(t)\;dt-\int_0^\infty f(t+i\varphi)\;dt=i
  \int_0^\varphi f(i t)\;dt\qquad \varphi\in(-\alpha,\beta).
\end{equation}
\end{proposition}
\begin{Proof} Given $\varphi\in(-\alpha,\beta)$, take $R\in\R$, with $R>0$, and consider de
rectangular contour $C$ of vertices $0$, $R$, $R+i\varphi$ and
$i\varphi$. For Cauchy theorem the integral of $f$ over $C$ must
be null; taking the limit for $R\to\infty$, the integral over the
right vertical side tends to zero and, parameterizing the
integrals over the other three sides, we have (\ref{ident:funol}).
\end{Proof}

\begin{definition}
Given a function $g(z)$, $M(g)(z)$ is its Mellin transform, that
is
\begin{equation}\label{def:mellintrasf}
  M(g)(z)=\int_0^\infty t^{z-1} g(t)\; dt
\end{equation}
for those values of $z$ for which the integral exists.
Furthermore, we indicate with $g_\varphi(z)$ its translated of a
value $i\varphi$ in the domain, that is
\begin{equation}\label{def:gphitrasl}
  g_\varphi(z)=g(z+i\varphi).
\end{equation}
\end{definition}

\begin{proposition}\label{teo:funolmel}Let $\alpha$ and $\beta$ be positive real numbers and
let $g(z)$ be an holomorphic function
in the strip $S=\left\{z\in\C\;:\;\Re(z)\geq 0, \\
\Im(z)\in(-\alpha,\beta)\right\}$ such that $g(z)=o(z^{-k})$, for
all $k\in\N$ if $\Re(z)\to\infty$. Then, for every integer $n\geq
0$
\begin{gather}\label{ident:funolmelpot}
  M(g)(n+1)-\sum_{k=0}^n
  \binom{n}{k}(i\varphi)^{n-k}M(g_\varphi)(k+1)=
  i\int_0^\varphi (i t)^n g(i t)\;dt \\
  \varphi\in(-\alpha,\beta).\nonumber
\end{gather}
Moreover, if $g(z)$ is not holomorphic in $z=0$, having a pole of
order $p$, the equality holds for every $n\geq p$.
\end{proposition}

\begin{Proof}
Let us consider $f(z)$ defined as
\begin{equation}
  f(z)=z^n g(z)
\end{equation}
where $n$ is chosen as specified in the above Proposition. Then
$f(z)$ satisfies the hypothesis of Proposition \ref{teo:funol}
and, in view of (\ref{def:gphitrasl}),
\begin{equation}
  f(z+i\varphi)=(z+i\varphi)^n g_{\varphi}(z).
\end{equation}
For the binomial theorem, we have
\begin{equation}
  f(z+i\varphi)= \sum_{k=0}^n \binom{n}{k}(i\varphi)^{n-k}z^k
  g_\varphi(z).
\end{equation}
Substituting in (\ref{ident:funol}), in view of
(\ref{def:mellintrasf}), noting that the Mellin transforms are
well defined, we obtain (\ref{ident:funolmelpot}).
\end{Proof}

\section{Recurrence relations for the Lerch $\Phi$ function}
\subsection{An important transform}\label{subsec:phimell} \ \\Given $a,b\in\C$
such that $|a|\leq1$ and $\Re(b)>0$, let us consider the function
of the complex variable $t$

\begin{equation}
  \eta(t)= \frac{e^{(1-b)t}}{e^{t}-a}.
\end{equation}
\\
Its Mellin transform (as defined in (\ref{def:mellintrasf})) is

\begin{equation}
  M(\eta)(s)=
  \int_0^\infty{\frac{t^{s-1}e^{(1-b)t}}{e^{t}-a}\,dt}.
\end{equation}
\\
We note that the integral exists if $\Re(s)>1$ when $a=1$ and if
$\Re(s)>0$ in the other cases. We can write

\begin{eqnarray}
M(\eta)(s)  & = &
\int_0^\infty{\frac{t^{s-1}e^{-bt}}{1-ae^{-t}}\,dt}\nonumber \\
            & = & \int_0^\infty{t^{s-1}e^{-bt}\sum_{k=0}^\infty
                  (ae^{-t})^k\,dt} \nonumber \\
            & = & \sum_{k=0}^\infty a^k \int_0^\infty
                  t^{s-1}e^{-(k+b)t}\;dt \nonumber \\
            & = & \sum_{k=0}^\infty
            \frac{a^k\Gamma(s)}{(k+b)^{s}} \nonumber
\end{eqnarray}
\\
so that we have (cf. for ex. \cite[pg. 121, eq. (4)]{SriChoi2001})

\begin{gather}\label{ident:PhiMel}
M(\eta)(s) = \Gamma(s)\Phi(a,s,b) \\
(|a|\leq1;\; \Re(b)>0;\; \Re(s)>1\; \mbox{when} \; a=1, \;
\Re(s)>0\; \mbox{otherwise})\nonumber.
\end{gather}

\subsection{Lerch $\Phi$ and Bernoulli polynomials series}
\
\\Referring to equation (\ref{ident:funolmelpot})
let us consider, for $b\in\C$ with $\Re(b)>0$, the function of the
complex variable $z$
\[
g(z)=\frac{e^{(1-b)z}}{e^z-1}\,,
\]
for which we have
\begin{align*}
g_\varphi(z)&=e^{-i b
\varphi}\frac{e^{(1-b)z}}{e^z-e^{-i\varphi}}\,.
\end{align*}
\\
Using (\ref{ident:PhiMel}), considering that
$\Phi(1,s,b)=\zeta(s,b)$, we easily find that, for $n\geq1$,
$M(g)(n+1)=n!\;\zeta(n+1,b)$ and, for $k\geq0$ and
$0<|\varphi|<2\pi$,
$M(g_\varphi)(k+1)=e^{-ib\varphi}k!\;\Phi(e^{-i\varphi},k+1,b)$.
\\
Furthermore, from  (\ref{def:PolBern}) we have
\begin{align*}
i\int_0^\varphi(it)^n g(it)\; dt     &=i\int_0^\varphi (it)^{n-1} \sum_{k=0}^\infty \frac{B_k(1-b)(it)^k}{k!}\\
                                     &=(i\varphi)^n\sum_{k=0}^\infty
                                     \frac{B_k(1-b)(i\varphi)^k}{k!\;(k+n)}\\
                                     &(n\geq1,\; |\varphi|<2\pi)
\end{align*}
So that (\ref{ident:funolmelpot}) becomes (multiplying both sides
by $(i\varphi)^{-n}$)
\begin{gather}
n!\;(i\varphi)^{-n}\zeta(n+1,b) - e^{-ib\varphi}\sum_{k=0}^n
\binom{n}{k}(i\varphi)^{-k}k!\;\Phi(e^{-i\varphi},k+1,b)=
\nonumber\\
\sum_{k=0}^\infty
\frac{B_k(1-b)(i\varphi)^k}{k!\;(k+n)}\label{ident:variamiameno}\\
(n\geq1,\;0<|\varphi|<2\pi, \Re(b)>0)\nonumber
\end{gather}

Some applications of this result, as well as its connection with
(\ref{EquSriv}), will be discussed in section \ref{applic}.

\subsection{Lerch $\Phi$ and Euler polynomials series}
\
\\Let us now consider, for $\Re(b)>0$, the function of the complex
variable $z$
\[
g(z)=\frac{e^{(1-b)z}}{e^z+1},
\]
\\
so that we have
\begin{align*}
g_\varphi(z)&=e^{-i b
\varphi}\frac{e^{(1-b)z}}{e^z+e^{-i\varphi}}.
\end{align*}
In view of (\ref{ident:PhiMel}), for $n\geq 0$,
$M(g)(n+1)=n!\;{\Phi(-1,n+1,b)}$ and, for $k\geq0$ and
$|\varphi|<\pi$,
$M(g_\varphi)(k+1)=e^{-ib\varphi}k!\;\Phi(-e^{-i\varphi},k+1,b)$.
Moreover, from (\ref{def:PolEul}) we have
\begin{align*}
i\int_0^\varphi(it)^n g(it)\; dt     &=i\int_0^\varphi \frac{(it)^{n}}{2} \sum_{k=0}^\infty \frac{E_k(1-b)(it)^k}{k!}\\
                                     &=\frac{(i\varphi)^{n+1}}{2}\sum_{k=0}^\infty \frac{E_k(1-b)(i\varphi)^k}{k!\;(k+n+1)}\\
                                     &(n\geq 0,\; |\varphi|<\pi).
\end{align*}
Substituting in (\ref{ident:funolmelpot}) and multiplying both
sides by $(i\varphi)^{-n}$ we obtain
\begin{gather}
n!\;(i\varphi)^{-n}{\Phi(-1,n+1,b)} - e^{-ib\varphi}\sum_{k=0}^n
\binom{n}{k}(i\varphi)^{-k}k!\;\Phi(-e^{-i\varphi},k+1,b)=
\nonumber\\\label{ident:variamiapiu} \frac{1}{2}\sum_{k=0}^\infty
\frac{E_k(1-b)(i\varphi)^{k+1}}{k!\;(k+n+1)}\\
(n\geq 0,\; 0<|\varphi|<\pi, \Re(b)>0).\nonumber
\end{gather}

\section{Applications of (\ref{ident:variamiameno})}\label{applic}

\subsection{Deduction of (\ref{EquSriv})}\
\\Suppose $b=1$ in (\ref{ident:variamiameno}); it can be easily shown
that $\zeta(n,1)=\zeta(n)$ and, if $k\geq0$ and
$0<|\varphi|<2\pi$, $\Phi(e^{-i \varphi},k+1,1)=e^{i \varphi}\;
\Li_{k+1}(e^{-i\varphi})$. So, equation (\ref{ident:variamiameno})
reduces to
\begin{gather}\label{ident:equasri01}
n!\;(i\varphi)^{-n}\zeta(n+1) - \sum_{k=0}^n
\binom{n}{k}(i\varphi)^{-k}k!\;\Li_{k+1}(e^{-i\varphi})=
\nonumber\\
\sum_{k=0}^\infty \frac{(i\varphi)^k B_k}{k!\;(k+n)}\\
(n\geq1,\;0<|\varphi|<2\pi)\nonumber
\end{gather}
\\
Moreover, it is well known (see \cite{handbook}) that $B_{1}=-1/2$,
$B_{2k+1}=0$ for $k>0$ and, if $\varphi \neq 0$, $\Li_1(e^{-i
\varphi})=-\log(1-e^{-i \varphi})$.

Thus, taking $\varphi=-2\pi/\omega$, equation
(\ref{ident:equasri01}) becomes
\begin{gather}\label{ident:equasri02}
n!\;\left(\frac{i\omega}{2\pi}\right)^n\zeta(n+1)+\log(1-e^{2\pi i
/\omega }) - \sum_{k=1}^n
\binom{n}{k}\left(\frac{i\omega}{2\pi}\right)^k
k!\;\Li_{k+1}(e^{2\pi i /\omega})=
\nonumber\\
\sum_{k=0}^\infty {\frac{(-1)^k(2\pi)^{2k}
B_{2k}}{(2k)!\;(2k+n)\omega^{2k}}}+\frac{\pi i}{(n+1)\omega}\\
(n\geq1,\;|\omega|>1).\nonumber
\end{gather}
Using (\ref{ident:zetanumbern}) to express $B_{2k}$ in terms of
$\zeta(2k)$, multiplying by $-1/2$ and rearranging terms, we
obtain
\begin{align}\label{EquSrivcor}
  \sum_{k=0}^\infty \frac{\zeta(2k)}{(2k+n)\omega^{2k}}& =\frac{\pi
                  i}{2(n+1)\omega}-\frac{1}{2}\log\left(1-e^{2\pi i
                  /\omega}\right)\\\nonumber
                & -\frac{n!}{2}\left(\frac{i\omega}{2
                \pi}\right)^n
                  \zeta(n+1)+\frac{1}{2}\sum_{k=1}^n \binom{n}{k}k!\left(\frac{i  \omega}{2
                  \pi}\right)^k\mbox{Li}_{k+1}\left( e^{2 \pi
                  i/\omega}\right)\\
                  & \mspace{100mu}(n\geq1,\;|\omega|>1)\nonumber
\end{align}
which is (\ref{EquSriv}) (and thus \cite[eq. 5.11]{SriGlaAdam2000}) with the correction of the index of the
PolyLogarithm.

As explained by the authors in \cite{SriGlaAdam2000}, equation
(\ref{EquSrivcor}) can be used to obtain several series
representation for the $\zeta$ for odd integer values of its
arguments. Some interesting examples, obtained for $\omega=2$, are (\cite[eq. 3.4 and 3.5]{SriGlaAdam2000})
\begin{equation}\label{zeta301}
\zeta(3)=\frac{2\pi^2}{9}\left(\log2 +2 \sum_{k=0}^\infty
\frac{\zeta(2k)}{2^{2k}(2k+3)}\right)\,.
\end{equation}
\begin{equation}\label{zeta302}
\zeta(3)=\frac{2\pi^2}{7}\left(\log2 +2 \sum_{k=0}^\infty
\frac{\zeta(2k)}{2^{2k}(2k+2)}\right)\,,
\end{equation}
Combining these two equations we have (\cite[eq. 2.19]{ChenSri1998} and \cite[eq. 3.32 with $n=1$]{Sri2000a})
\begin{equation}\label{zeta303}
\zeta(3)=-2\pi^2\sum_{k=0}^\infty
\frac{\zeta(2k)}{2^{2k}(2k+2)(2k+3)}
\end{equation}
and
\begin{equation}\label{log201}
\log2=-\sum_{k=0}^\infty
\frac{\zeta(2k)(4k+13)}{2^{2k}(2k+2)(2k+3)}\,.
\end{equation}
On the other hand, when $\omega=4$, more interesting results can
be obtained. It is easy to see that
\begin{equation}\label{ident:polizetabeta}
  \Li_k(i)=\left(\frac{1-2^{k-1}}{2^{2k-1}}\right)\zeta(k)+i
  \beta(k)\qquad (k\geq2).
\end{equation}
Substituting in (\ref{EquSrivcor}), taking the real part of both
sides (paying attention on the parity of $n$) it is possible to
obtain the following results (calculations are omitted for
brevity)
\begin{multline}\label{ident:bordel1}
\sum_{k=0}^n \binom{2n+1}{2k+1}\frac{2^{2k+1}(-1)^k
(2k+1)!}{\pi^{2k+1}}\;\beta(2k+2)\\
+\sum_{k=1}^n \binom{2n+1}{2k}\frac{2^{2k}(2^{2k}-1)(-1)^k
(2k)!}{2^{4k+1}\pi^{2k}}\;\zeta(2k+1)=\\ -\frac{\log
2 }{2}-2\sum_{k=0}^\infty \frac{\zeta(2k)}{2^{4k}(2k+2n+1)}\qquad
(n\geq0)
\end{multline}
and
\begin{multline}\label{ident:bordel2}
(-1)^n\left(\frac{2}{\pi}\right)^{2n}(2n)!\;\zeta(2n+1) \\
+\sum_{k=1}^n \binom{2n}{2k-1}\frac{\pi^{1-2k}(-1)^{1-k}
(2k-1)!}{2^{1-2k}}\;\beta(2k)\\
+\sum_{k=1}^n \binom{2n}{2k}\frac{2^{2k}(2^{2k}-1)(-1)^k
(2k)!}{2^{4k+1}\pi^{2k}}\;\zeta(2k+1)=\\ -\frac{\log 2
}{2}-2\sum_{k=0}^\infty \frac{\zeta(2k)}{2^{4k}(2k+2n)}\qquad
(n\geq1),
\end{multline}
which are equivalent to \cite[eq. 3.26 and 3.27]{SriGlaAdam2000} but where we have isolated the $\beta$ function terms (which corresponds to express the Clausen's functions of \cite[eq. 2.18 and 2.19]{SriGlaAdam2000}, when $\omega=4$, in terms of the $\beta$ function instead of using \cite[eq. 3.7 and 3.16]{SriGlaAdam2000}).

If, in (\ref{ident:bordel1}), $n=0$, we have (\cite[eq. 2.23]
{SriGlaAdam2000})
\begin{equation}\label{L24k}
 \beta(2)=\mbox{Catalan}=-\frac{\pi}{4}\left(\log(2)+4\sum_{k=0}^\infty
\frac{\zeta(2k)}{2^{4k}(2k+1)}\right)
\end{equation}
and substituting this representation in (\ref{ident:bordel2}) with
$n=1$
\begin{equation}\label{Z34k}
\zeta(3)=-\frac{2\pi^2}{35}\left(\log(2)+ 4\sum_{k=0}^\infty
\frac{\zeta(2k)(2k+3)}{2^{4k}(2k+1)(2k+2)}\right).
\end{equation}
The formula (\ref{Z34k}) is essentially the same as a known result \cite[p. 192, Equation (3.21)]{ChenSri1998}. Applying alternatively (\ref{ident:bordel1}) and
(\ref{ident:bordel2}) with increasing values of $n$ and using the
representations found at every step we have
\begin{equation}\label{L44k}
\beta(4)=-\frac{\pi^3}{10080}\left(183\log(2)+12\sum_{k=0}^\infty
\frac{\zeta(2k)(244k^2+732k+479)}{2^{4k}(2k+1)(2k+2)(2k+3)}\right),
\end{equation}
\bigskip
\begin{equation}\label{Z54k}
\zeta(5)=-\frac{\pi^4}{166005}\left(942\log(2)+48\sum_{k=0}^\infty
\frac{\zeta(2k)(628k^3+3140k^2+5111k+2581)}{2^{4k}(2k+1)(2k+2)(2k+3)(2k+4)}\right)
\end{equation}
and so on.

\subsection{Another application: case b=1/2}\
\\We present here another possible use of equation
(\ref{ident:variamiameno}) when $b=1/2$ and $\varphi=\pi$.
\\It is easy to verify that $\zeta(n,1/2)=(2^n-1)\zeta(n)$
and, if $k\geq 0$, $\Phi(e^{-i \pi},k+1,1/2)=2^{k+1}\beta(k+1)$.
\\
It is also easy to demonstrate that $B_{k}(1/2)  =
(2^{1-k}-1)B_{k}$; thus equation (\ref{ident:variamiameno})
reduces to
\begin{gather}
n!\;(i \pi)^{-n}(2^{n+1}-1)\zeta(n+1) + i \sum_{k=0}^n
\binom{n}{k}(i \pi)^{-k}k!\;2^{k+1}\beta(k+1)=
\nonumber\\
\sum_{k=0}^\infty
\frac{(2^{1-2k}-1)B_{2k}\;(-1)^k\pi^{2k}}{(2k)!\;(2k+n)}\\
(n\geq1)\nonumber
\end{gather}
and using (\ref{ident:zetanumbern}) we have
\begin{gather}
n!\;(i \pi)^{-n}(2^{n+1}-1)\zeta(n+1) + i \sum_{k=0}^n
\binom{n}{k}(i \pi)^{-k}k!\;2^{k+1}\beta(k+1)=
\nonumber\\
\sum_{k=0}^\infty
\frac{(2^{2k-1}-1)\zeta(2k)}{2^{4k-2}(2k+n)}\label{ident:convarmiameno}\\
(n\geq1).\nonumber
\end{gather}

This equation allows us to obtain some relations; taking the
imaginary part of both side we have the following two formulas
\begin{equation}\label{ident:betarecurel}
  \sum_{k=0}^n  \frac{2^{2k}(-1)^k}{(2n-2k)!\;\pi^{2k}}\;
  \beta(2k+1)=0\qquad (n\geq1)
\end{equation}
%\bigskip
\begin{gather}\label{ident:zetacombeta}
  \zeta(2n+2)=\frac{(-1)^n \pi^{2n+1}}{(2^{2n+2}-1)}\sum_{k=0}^n
  \frac{2^{2k+1}(-1)^k }{(2n-2k+1)!\;\pi^{2k}}\;\beta(2k+1)\\
  (n\geq0)\nonumber
\end{gather}
which, in view of (\ref{ident:zetanumbern}) and
(\ref{ident:betanumeul}), show to be equivalent to (cf. for ex. \cite[pg. 64 eq (48)]{SriChoi2001})
\begin{equation}\label{ident:numeulrecurel}
\sum_{k=0}^{n}\binom{2n}{2k}E_{2k}=0\qquad(n\geq1)
\end{equation}
and
\begin{equation}\label{ident:numbercombnumeul}
\sum_{k=0}^{n-1}\binom{2n-1}{2k}E_{2k}=\frac{2^{2n}(2^{2n}-1)}{2n}B_{2n}\qquad(n\geq1)
\end{equation}
respectively (or, \emph{viceversa}, known (\ref{ident:numeulrecurel}) and (\ref{ident:numbercombnumeul}), (\ref{ident:betarecurel}) and (\ref{ident:zetacombeta}) are equivalent to (\ref{ident:zetanumbern}) and
(\ref{ident:betanumeul})).

Now, taking the real part of both sides of
(\ref{ident:convarmiameno}) we have
\begin{gather}\label{ident:betaserepr}
  \sum_{k=0}^n \binom{2n+1}{2k+1}\frac{2^{2k+1}(2k+1)!}{\pi^{2k+1}}\;(-1)^k \beta(2k+2)=\sum_{k=0}^\infty
  \frac{(2^{2k-1}-1)\zeta(2k)}{2^{4k-1}(2k+2n+1)}\\
  (n\geq0)\nonumber
\end{gather}
\begin{multline}\label{ident:zetaserepr}
  \zeta(2n+1)=\frac{(-1)^n \pi^{2n}}{(2n)!(2^{2n+1}-1)} \left( \sum_{k=0}^\infty \frac{(2^{2k-1}-1)\zeta(2k)}{2^{4k-2}(2k+2n)} \right. \\ \left.-  \sum_{k=0}^{n-1} \binom{2n}{2k+1}\frac{2^{2k+2}(2k+1)!}{\pi^{2k+1}}\;(-1)^k
  \beta(2k+2)\right)\qquad(n\geq1)
\end{multline}

From equation (\ref{ident:betaserepr}), some interesting series
representations for $\beta(2m)$, $m\in \N$ can be derived. When
$n=0$, for example, we have
\begin{equation}\label{beta(2)}
  \beta(2)=\pi \sum_{k=0}^{\infty}\frac{(2^{2k-1}-1)\zeta(2k)}{2^{4k}(2k+1)}\,,
\end{equation}
which can also be obtained by combining $\ref{L24k}$ and the known (see \cite{Sri_Taiwan2000}, eq. 4.11, pg. 586) sum (\ref{sumzetalog2}), while for $n=1$ we have
\begin{equation}\label{beta(4)}
  \beta(4)=\frac{ \pi^3}{6} \sum_{k=0}^{\infty}\frac{(2^{2k-1}-1)(k+2) \zeta(2k)}
  {2^{4k}(2k+1)(2k+3)}\,.
\end{equation}
Using (\ref{ident:zetaserepr}), similar series representations for
the values $\zeta(2m+1)$ $m\in\N$ can be obtained, for example
\begin{equation}\label{Z(3)}
  \zeta(3)=\frac{2 \pi^2}{7} \sum_{k=0}^{\infty}\frac{(2^{2k-1}-1)(2k+3) \zeta(2k)}
  {2^{4k}(2k+1)(2k+2)}\,,
\end{equation}
\begin{equation}\label{Z(5)}
  \zeta(5)=\frac{ \pi^4}{186} \sum_{k=0}^{\infty}\frac{(2^{2k-1}-1)(20k^2+80k+83) \zeta(2k)}
  {2^{4k}(2k+1)(2k+3)(2k+4)}\,.
\end{equation}
In view of (\ref{beta(2)}), the formula (\ref{Z(3)}) is equivalent to a known result \cite[p. 191, Equation (3.13)]{ChenSri1998}.

\section{Applications of (\ref{ident:variamiapiu})}
\subsection{Series of $\beta$: companion of (\ref{EquSrivcor})}\
\\
From equation (\ref{ident:variamiapiu}), if $b=1/2$, a companion
of equation (\ref{EquSrivcor}) can be derived, in which the rule
of the $\zeta$ function is played by the $\beta$ function.
\\It is easy to see that $\Phi(-1,n+1,1/2)=2^{n+1}\beta(n+1)$.
Now, taking $\varphi=-\pi/\omega$, using (\ref{ident:betanumeul})
and remembering that $E_k=2^k E_k(1/2)$, we can rewrite the right
hand side of (\ref{ident:variamiapiu}) as
\begin{equation}
  -2i \sum_{k=0}^\infty \frac{\beta(2k+1)}{\omega^{2k+1}(2k+n+1)}\,,
\end{equation}
where we have used the fact that $E_{2k+1}=0$.
\\Thus, multiplying
both sides by $-1/2$, equation (\ref{ident:variamiapiu}) can be
rewritten as
\begin{multline}\label{Equsrivvar}
i \sum_{k=0}^\infty \frac{\beta(2k+1)}{\omega^{2k+1}(2k+n+1)}=
  -\, n!\; \left(\frac{2\omega i}{\pi}\right)^n
  \beta(n+1)+ \\ \frac{1}{2}e^{i\pi/(2\omega)}\sum_{k=0}^n
\binom{n}{k}\left(\frac{i\omega}{\pi}\right)^{k}
k!\;\Phi(-e^{i\pi/\omega},k+1,1/2)\mspace{-150mu}\\(n\geq1,\;|\omega|<1).\mspace{200mu}
\end{multline}

\subsection{Case $\omega=-2$ and Dirichlet $L$-series}\
\\As an example of application of equation (\ref{Equsrivvar}), let
us set $\omega=-2$ so that, multiplying both sides for $-2$, the
left hand side becomes
\begin{equation}
  i \sum_{k=0}^\infty \frac{\beta(2k+1)}{2^{2k}(2k+n+1)}.
\end{equation}
With some further calculations, it is possible to demonstrate that
for every $s>1$
\begin{equation}
\Phi(i,s,1/2)=2^{-2s}(\zeta(s,1/8)-\zeta(s,5/8)+
i\zeta(s,3/8)-i\zeta(s,7/8))
\end{equation}
so that we have
\begin{equation}\label{ident:phiLseries}
  e^{-i \pi/4}\Phi(i,s,1/2)=2^{s-1/2}\left(L(s,\chi_1)-i L(s,\chi_2)
  \right)\,,
\end{equation}
where $\chi_1$ and $\chi_2$ are characters on $\Z/8$ satisfying
\begin{equation}
\begin{array}{ccc}
  \chi_1(n)=\left\{
  \begin{array}{rl}
    1 & \text{if}\quad n=1,3 \\
    -1 & \text{if}\quad n=5,7
  \end{array}\right.
                       & \qquad   & \chi_2(n)= \left\{
  \begin{array}{rl}
    1 & \text{if}\quad n=1,7 \\
    -1 & \text{if}\quad n=3,5.
  \end{array}\right.
\end{array}
\end{equation}
Considering the convergence of the $L$-series, we can say that
equation (\ref{ident:phiLseries}) holds for every integer $s=k$
with $k\geq1$. As a result, equation (\ref{Equsrivvar}) becomes
\begin{multline}
  n!\;\left(i \pi\right)^{-n}2^{2n+1}\beta(n+1)-\\ \sum_{k=0}^n
  \binom{n}{k}k!\; (i \pi)^{-k}2^{2k+1/2}\left(L(k+1,\chi_1)-i L(k+1,\chi_2)
  \right)= \\
  i \sum_{k=0}^\infty
  \frac{\beta(2k+1)}{2^{2k}(2k+n+1)}\qquad\qquad(n\geq0).
  \label{ident:betaLbeta}
\end{multline}
Taking the imaginary part of both sides we have the following two
relations (see \cite{SriTsu2001} for interesting more general, but not equivalent, recursions and series representation for Dirichlet $L$-series)
\begin{multline}\label{ident:bordelL01}
\sum_{k=0}^{n-1}\binom{2n}{2k+1}\frac{(-1)^k
2^{4k+5/2}(2k+1)!}{\pi^{2k+1}}\;L(2k+2,\chi_1)+\\
\sum_{k=0}^{n}\binom{2n}{2k}\frac{(-1)^k
2^{4k+1/2}(2k)!}{\pi^{2k}}\;L(2k+1,\chi_2)=\\
\sum_{k=0}^\infty
\frac{\beta(2k+1)}{2^{2k}(2k+2n+1)}\quad(n\geq0),
\end{multline}
\begin{multline}\label{ident:bordelL02}
\frac{(-1)^{n-1}2^{4n+3}(2n+1)!}{\pi^{2n+1}}\beta(2n+2)+\\
\sum_{k=0}^{n}\binom{2n+1}{2k+1}\frac{(-1)^k
2^{4k+5/2}(2k+1)!}{\pi^{2k+1}}\;L(2k+2,\chi_1)+\\
\sum_{k=0}^{n}\binom{2n+1}{2k}\frac{(-1)^k
2^{4k+1/2}(2k)!}{\pi^{2k}}\;L(2k+1,\chi_2)=\\
\sum_{k=0}^\infty \frac{\beta(2k+1)}{2^{2k}(2k+2n+2)}\quad(n\geq0)
\end{multline}
If, for example, on equation (\ref{ident:bordelL01}), $n=0$, we
have
\begin{equation}\label{ident:caso0L}
  \sqrt{2}\,L(1,\chi_2)=\sum_{k=0}^\infty
  \frac{\beta(2k+1)}{2^{2k}(2k+1)}\,.
\end{equation}
On the other hand one may verify that
\begin{equation}\nonumber
L(1,\chi_2)=\frac{\sqrt 2}{2} \sum_{n=1}^\infty
\frac{1-(-1)^n}{n}\cos(n\pi/4)
\end{equation}
so that, writing the cosines in exponential form, we can sum the
series to obtain $L(1,\chi_2)=\log(1+\sqrt 2)/\sqrt 2$. This gives
us the result
\begin{equation}\label{ident:caso0Lval}
  \sum_{k=0}^\infty
  \frac{\beta(2k+1)}{2^{2k}(2k+1)}=\log(1+\sqrt 2)\,,
\end{equation}
which is an interesting (presumably new) counterpart of the known sum
(\cite{Sri_Taiwan2000}, eq. 4.11, pg. 586)
\begin{equation}\label{sumzetalog2}
  \sum_{k=0}^\infty
  \frac{\zeta(2k)}{2^{2k}(2k+1)}=-\frac{1}{2}\log(2)\,.
\end{equation}

These examples show some possible uses of equations
(\ref{ident:variamiameno}) and (\ref{ident:variamiapiu}), and it
is clear that other chooses of the parameters $b$ and $\varphi$
will give more complicated, but maybe interesting, equalities and
recurrence relations.

\section{Reflection properties}
It is well known that Bernoulli and Euler polynomials satisfy the
property
\begin{equation}
B_k(1-x)=(-1)^k B_k(x)
\end{equation}
and
\begin{equation}
E_k(1-x)=(-1)^k E_k(x).
\end{equation}
So, called $S_1(n,\varphi,b)$ the left hand side of
(\ref{ident:variamiameno}), we have
\begin{gather}
  S_1(n,\varphi,b)=\overline{S_1(n,\varphi,1-b)}=S_1(n,-\varphi,1-b)\\
(n\geq1,\;0<|\varphi|<2\pi,\;0<\Re(b)<1)\nonumber
\end{gather}
and, called $S_2(n,\varphi,b)$ the left hand side of
(\ref{ident:variamiapiu}), we have
\begin{gather}
  S_2(n,\varphi,b)=-\overline{S_2(n,\varphi,1-b)}=-S_2(n,-\varphi,1-b)\\
(n\geq1,\;0<|\varphi|<2\pi,\;0<\Re(b)<1)\nonumber
\end{gather}
We have found a relation between the values $\Phi(e^{i
\varphi},k,b)$ and $\Phi(e^{i \varphi},k,1-b)$, $k=1\ldots n$,
that does not include any infinite sum.

\subsection*{Acknowledgements.\\}
I thank very much Gabriella Rossi, Andrea Marson and Alessandro
Languasco for helping me.


\begin{thebibliography}{99}

\bibitem{handbook} Abramowitz,~M.~and I.~Stegun:
{\it Handbook of Mathematical Functions.} Dover, N.~Y., 1972.

\bibitem{AdamSri1998} Adamchik,~V.~S.~and H. M. Srivastava~:
{\it Some series of the Zeta and related functions.}
Analysis~18(1998), 131-144.

\bibitem{ChenSri1998} Chen,~M.~-P.~and H. M. Srivastava~:
{\it Some families of series representations for the Riemann $\zeta(3)$.}
Resultate Math.~33(1998), 179-197.

\bibitem{Ramasw} Ramaswami,~V.:
{\it Notes on Riemann's $\zeta$-function.} J. London Math
Soc.~9(1943), 165-169.

\bibitem{Sri2000a} Srivastava,~H.~M.:
{\it Some simple algorithms for the evaluations and representations of the Riemann Zeta function at positive integer arguments.} J. Math. Anal. Appl.~246(2000),
331-351.

\bibitem{Sri_Taiwan2000} Srivastava,~H.~M.:
{\it Some families of rapidly convergent series representations
for the zeta function.} Taiwanese J. of Math.~4(2000),
569-598.

\bibitem{SriChoi2001} Srivastava,~H.~M.~and J.~Choi:
{\it Series Associated with the Zeta and Related Functions.} Kluwer Academic Publishers, Dordrecht, Boston, and London, 2001.

\bibitem{SriGlaAdam2000} Srivastava,~H.~M.~, Glasser,~M.~L.~and V.~S.~Adamchik:
{\it Some definite integrals associated with the riemann zeta
function.} Z.~Anal. Anwendungen~19(2000), 831-846.

\bibitem{SriTsu2001} Srivastava,~H.~M.~and H.~Tsumura:
{\it Certain classes of rapidly convergent series representations for $L(2n,\chi)$ and $L(2n+1,\chi)$.} Acta Arith.~100(2001), 195-201.
\end{thebibliography}
\end{document}